\documentclass[12pt]{amsart}
\usepackage{amsmath,amssymb,amsbsy,amsfonts,amsthm,latexsym,
                        amsopn,amstext,amsxtra,euscript,amscd,mathrsfs,color,bm}
                        
\usepackage{float} 
\usepackage[english]{babel}
\usepackage{mathtools}
\usepackage[textsize=tiny]{todonotes}
\usepackage{url}
\usepackage[colorlinks,linkcolor=blue,anchorcolor=blue,citecolor=blue]{hyperref}
\usepackage{nicefrac}
\usepackage{enumitem}
\usepackage{comment}
\usepackage{orcidlink}

\usepackage[a4paper,  left=3.6cm, right=3.6cm, top=3.5 cm, bottom=3.5cm]{geometry}

\usepackage[norefs,nocites]{refcheck}
     
  \usepackage{soul}

\newcommand{\rvline}{\hspace*{-\arraycolsep}\vline\hspace*{-\arraycolsep}}

\begin{document}
%
    
\newtheorem{theorem}{Theorem}
\newtheorem{lemma}[theorem]{Lemma}
\newtheorem{example}[theorem]{Example}
\newtheorem{algol}{Algorithm}
\newtheorem{corollary}[theorem]{Corollary}
\newtheorem{prop}[theorem]{Proposition}
\newtheorem{proposition}[theorem]{Proposition}
\newtheorem{problem}[theorem]{Problem}
\newtheorem{conj}[theorem]{Conjecture}

\theoremstyle{remark}
\newtheorem{definition}[theorem]{Definition}
\newtheorem{question}[theorem]{Question}
\newtheorem{remark}[theorem]{Remark}
\newtheorem*{acknowledgement}{Acknowledgements}

\newtheorem*{Thm*}{Theorem}
\newtheorem{Thm}{Theorem}[section]
\renewcommand*{\theThm}{\Alph{Thm}}

\numberwithin{equation}{section}
\numberwithin{theorem}{section}
\numberwithin{table}{section}
\numberwithin{figure}{section}

\allowdisplaybreaks

\definecolor{olive}{rgb}{0.3, 0.4, .1}
\definecolor{dgreen}{rgb}{0.,0.5,0.}

\def\cc#1{\textcolor{red}{#1}} 

\definecolor{dgreen}{rgb}{0.,0.6,0.}
\def\tgreen#1{\begin{color}{dgreen}{\it{#1}}\end{color}}
\def\tblue#1{\begin{color}{blue}{\it{#1}}\end{color}}
\def\tred#1{\begin{color}{red}#1\end{color}}
\def\tmagenta#1{\begin{color}{magenta}{\it{#1}}\end{color}}
\def\tNavyBlue#1{\begin{color}{NavyBlue}{\it{#1}}\end{color}}
\def\tMaroon#1{\begin{color}{Maroon}{\it{#1}}\end{color}}

%


 \def\mand{\qquad\mbox{and}\qquad}

\def\cA{{\mathcal A}}
\def\cB{{\mathcal B}}
\def\cC{{\mathcal C}}
\def\cD{{\mathcal D}}
\def\cH{{\mathcal H}}
\def\cI{{\mathcal I}}
\def\cJ{{\mathcal J}}
\def\cK{{\mathcal K}}
\def\cL{{\mathcal L}}
\def\cM{{\mathcal M}}
\def\cN{{\mathcal N}}
\def\cO{{\mathcal O}}
\def\cP{{\mathcal P}}
\def\cQ{{\mathcal Q}}
\def\cR{{\mathcal R}}
\def\cS{{\mathcal S}}
\def\cT{{\mathcal T}}
\def\cU{{\mathcal U}}
\def\cV{{\mathcal V}}
\def\cW{{\mathcal W}}
\def\cX{{\mathcal X}}
\def\cY{{\mathcal Y}}
\def\cZ{{\mathcal Z}}

\def\C{\mathbb{C}}
\def\F{\mathbb{F}}
\def\K{\mathbb{K}}
\def\Z{\mathbb{Z}}
\def\R{\mathbb{R}}
\def\Q{\mathbb{Q}}
\def\N{\mathbb{N}}
\def\M{\mathrm{M}}
\def\L{\mathbb{L}}
\def\M{{\normalfont\textsf{M}}} 
\def\U{\mathbb{U}}
\def\P{\mathbb{P}}
\def\A{\mathbb{A}}
\def\fp{\mathfrak{p}}
\def\fq{\mathfrak{q}}
\def\n{\mathfrak{n}}
\def\X{\mathcal{X}}
\def\x{\textrm{\bf x}}
\def\w{\textrm{\bf w}}
\def\ovQ{\overline{\Q}}
\def \Kab{\K^{\mathrm{ab}}}
\def \Qab{\Q^{\mathrm{ab}}}
\def \Qtr{\Q^{\mathrm{tr}}}
\def \Kc{\K^{\mathrm{c}}}
\def \Qc{\Q^{\mathrm{c}}}
\def\ZK{\Z_\K}
\def\ZKS{\Z_{\K,\cS}}
\def\ZKSf{\Z_{\K,\cS_{f}}}
\def\RSf{R_{\cS_{f}}}
\def\RTf{R_{\cT_{f}}}

\def\S{\mathcal{S}}
\def\vec#1{\mathbf{#1}}
\def\ov#1{{\overline{#1}}}
\def\sign{{\operatorname{sign}}}
\def\Gm{\G_{\textup{m}}}
\def\fA{{\mathfrak A}}
\def\fB{{\mathfrak B}}

\def \GL{\mathrm{GL}}
\def \Mat{\mathrm{Mat}}

\def\house#1{{%
    \setbox0=\hbox{$#1$}
    \vrule height \dimexpr\ht0+1.4pt width .5pt depth \dp0\relax
    \vrule height \dimexpr\ht0+1.4pt width \dimexpr\wd0+2pt depth \dimexpr-\ht0-1pt\relax
    \llap{$#1$\kern1pt}
    \vrule height \dimexpr\ht0+1.4pt width .5pt depth \dp0\relax}}


\newenvironment{notation}[0]{%
  \begin{list}%
    {}%
    {\setlength{\itemindent}{0pt}
     \setlength{\labelwidth}{1\parindent}
     \setlength{\labelsep}{\parindent}
     \setlength{\leftmargin}{2\parindent}
     \setlength{\itemsep}{0pt}
     }%
   }%
  {\end{list}}

\newenvironment{parts}[0]{%
  \begin{list}{}%
    {\setlength{\itemindent}{0pt}
     \setlength{\labelwidth}{1.5\parindent}
     \setlength{\labelsep}{.5\parindent}
     \setlength{\leftmargin}{2\parindent}
     \setlength{\itemsep}{0pt}
     }%
   }%
  {\end{list}}
\newcommand{\Part}[1]{\item[\upshape#1]}

\def\Case#1#2{%
\smallskip\paragraph{\textbf{\boldmath Case #1: #2.}}\hfil\break\ignorespaces}

\def\Subcase#1#2{%
\smallskip\paragraph{\textit{\boldmath Subcase #1: #2.}}\hfil\break\ignorespaces}

\renewcommand{\a}{\alpha}
\renewcommand{\b}{\beta}
\newcommand{\g}{\gamma}
\renewcommand{\d}{\delta}
\newcommand{\e}{\epsilon}
\newcommand{\f}{\varphi}
\newcommand{\fhat}{\hat\varphi}
\newcommand{\bfphi}{{\boldsymbol{\f}}}
\renewcommand{\l}{\lambda}
\renewcommand{\k}{\kappa}
\newcommand{\lhat}{\hat\lambda}
\newcommand{\bfmu}{{\boldsymbol{\mu}}}
\renewcommand{\o}{\omega}
\renewcommand{\r}{\rho}
\newcommand{\rbar}{{\bar\rho}}
\newcommand{\s}{\sigma}
\newcommand{\sbar}{{\bar\sigma}}
\renewcommand{\t}{\tau}
\newcommand{\z}{\zeta}


\newcommand{\ga}{{\mathfrak{a}}}
\newcommand{\gb}{{\mathfrak{b}}}
\newcommand{\gn}{{\mathfrak{n}}}
\newcommand{\gp}{{\mathfrak{p}}}
\newcommand{\gP}{{\mathfrak{P}}}
\newcommand{\gq}{{\mathfrak{q}}}
\newcommand{\h}{{\mathfrak{h}}}
\newcommand{\Abar}{{\bar A}}
\newcommand{\Ebar}{{\bar E}}
\newcommand{\kbar}{{\bar k}}
\newcommand{\Kbar}{{\bar K}}
\newcommand{\Pbar}{{\bar P}}
\newcommand{\Sbar}{{\bar S}}
\newcommand{\Tbar}{{\bar T}}
\newcommand{\gbar}{{\bar\gamma}}
\newcommand{\lbar}{{\bar\lambda}}
\newcommand{\ybar}{{\bar y}}
\newcommand{\phibar}{{\bar\f}}

\newcommand{\Acal}{{\mathcal A}}
\newcommand{\Bcal}{{\mathcal B}}
\newcommand{\Ccal}{{\mathcal C}}
\newcommand{\Dcal}{{\mathcal D}}
\newcommand{\Ecal}{{\mathcal E}}
\newcommand{\Fcal}{{\mathcal F}}
\newcommand{\Gcal}{{\mathcal G}}
\newcommand{\Hcal}{{\mathcal H}}
\newcommand{\Ical}{{\mathcal I}}
\newcommand{\Jcal}{{\mathcal J}}
\newcommand{\Kcal}{{\mathcal K}}
\newcommand{\Lcal}{{\mathcal L}}
\newcommand{\Mcal}{{\mathcal M}}
\newcommand{\Ncal}{{\mathcal N}}
\newcommand{\Ocal}{{\mathcal O}}
\newcommand{\Pcal}{{\mathcal P}}
\newcommand{\Qcal}{{\mathcal Q}}
\newcommand{\Rcal}{{\mathcal R}}
\newcommand{\Scal}{{\mathcal S}}
\newcommand{\Tcal}{{\mathcal T}}
\newcommand{\Ucal}{{\mathcal U}}
\newcommand{\Vcal}{{\mathcal V}}
\newcommand{\Wcal}{{\mathcal W}}
\newcommand{\Xcal}{{\mathcal X}}
\newcommand{\Ycal}{{\mathcal Y}}
\newcommand{\Zcal}{{\mathcal Z}}

\renewcommand{\AA}{\mathbb{A}}
\newcommand{\BB}{\mathbb{B}}
\newcommand{\CC}{\mathbb{C}}
\newcommand{\FF}{\mathbb{F}}
\newcommand{\G}{\mathbb{G}}
\newcommand{\KK}{\mathbb{K}}
\newcommand{\NN}{\mathbb{N}}
\newcommand{\PP}{\mathbb{P}}
\newcommand{\QQ}{\mathbb{Q}}
\newcommand{\RR}{\mathbb{R}}
\newcommand{\ZZ}{\mathbb{Z}}

\newcommand{\bfa}{{\boldsymbol a}}
\newcommand{\bfb}{{\boldsymbol b}}
\newcommand{\bfc}{{\boldsymbol c}}
\newcommand{\bfd}{{\boldsymbol d}}
\newcommand{\bfe}{{\boldsymbol e}}
\newcommand{\bff}{{\boldsymbol f}}
\newcommand{\bfg}{{\boldsymbol g}}
\newcommand{\bfi}{{\boldsymbol i}}
\newcommand{\bfj}{{\boldsymbol j}}
\newcommand{\bfk}{{\boldsymbol k}}
\newcommand{\bfm}{{\boldsymbol m}}
\newcommand{\bfp}{{\boldsymbol p}}
\newcommand{\bfr}{{\boldsymbol r}}
\newcommand{\bfs}{{\boldsymbol s}}
\newcommand{\bft}{{\boldsymbol t}}
\newcommand{\bfu}{{\boldsymbol u}}
\newcommand{\bfv}{{\boldsymbol v}}
\newcommand{\bfw}{{\boldsymbol w}}
\newcommand{\bfx}{{\boldsymbol x}}
\newcommand{\bfy}{{\boldsymbol y}}
\newcommand{\bfz}{{\boldsymbol z}}
\newcommand{\bfA}{{\boldsymbol A}}
\newcommand{\bfF}{{\boldsymbol F}}
\newcommand{\bfB}{{\boldsymbol B}}
\newcommand{\bfD}{{\boldsymbol D}}
\newcommand{\bfG}{{\boldsymbol G}}
\newcommand{\bfI}{{\boldsymbol I}}
\newcommand{\bfM}{{\boldsymbol M}}
\newcommand{\bfP}{{\boldsymbol P}}
\newcommand{\bfX}{{\boldsymbol X}}
\newcommand{\bfY}{{\boldsymbol Y}}
\newcommand{\bfzero}{{\boldsymbol{0}}}
\newcommand{\bfone}{{\boldsymbol{1}}}

\newcommand{\aff}{{\textup{aff}}}
\newcommand{\Aut}{\operatorname{Aut}}
\newcommand{\Berk}{{\textup{Berk}}}
\newcommand{\Birat}{\operatorname{Birat}}
\newcommand{\characteristic}{\operatorname{char}}
\newcommand{\codim}{\operatorname{codim}}
\newcommand{\Crit}{\operatorname{Crit}}
\newcommand{\critwt}{\operatorname{critwt}} 
\newcommand{\cond}{\operatorname{cond}}
\newcommand{\Cycle}{\operatorname{Cycles}}
\newcommand{\diag}{\operatorname{diag}}
\newcommand{\Disc}{\operatorname{Disc}}
\newcommand{\Div}{\operatorname{Div}}
\newcommand{\Dom}{\operatorname{Dom}}
\newcommand{\End}{\operatorname{End}}
\newcommand{\ExtOrbit}{\mathcal{EO}} 
\newcommand{\Fbar}{{\bar{F}}}
\newcommand{\Fix}{\operatorname{Fix}}
\newcommand{\FOD}{\operatorname{FOD}}
\newcommand{\FOM}{\operatorname{FOM}}
\newcommand{\Gal}{\operatorname{Gal}}
\newcommand{\genus}{\operatorname{genus}}
\newcommand{\GITQuot}{/\!/}
\newcommand{\GR}{\operatorname{\mathcal{G\!R}}}
\newcommand{\Hom}{\operatorname{Hom}}
\newcommand{\Index}{\operatorname{Index}}
\newcommand{\Image}{\operatorname{Image}}
\newcommand{\Isom}{\operatorname{Isom}}
\newcommand{\hhat}{{\hat h}}
\newcommand{\Ker}{{\operatorname{ker}}}
\newcommand{\Ksep}{K^{\textup{sep}}}  
\newcommand{\lcm}{{\operatorname{lcm}}}
\newcommand{\LCM}{{\operatorname{LCM}}}
\newcommand{\Lift}{\operatorname{Lift}}
\newcommand{\limstar}{\lim\nolimits^*}
\newcommand{\limstarn}{\lim_{\hidewidth n\to\infty\hidewidth}{\!}^*{\,}}
\newcommand{\llog}{\log\log}
\newcommand{\logplus}{\log^{\scriptscriptstyle+}}
\newcommand{\maxplus}{\operatornamewithlimits{\textup{max}^{\scriptscriptstyle+}}}
\newcommand{\MOD}[1]{~(\textup{mod}~#1)}
\newcommand{\Mor}{\operatorname{Mor}}
\newcommand{\Moduli}{\mathcal{M}}
\newcommand{\Norm}{{\operatorname{\mathsf{N}}}}
\newcommand{\notdivide}{\nmid}
\newcommand{\normalsubgroup}{\triangleleft}
\newcommand{\NS}{\operatorname{NS}}
\newcommand{\onto}{\twoheadrightarrow}
\newcommand{\ord}{\operatorname{ord}}
\newcommand{\Orbit}{\mathcal{O}}
\newcommand{\Per}{\operatorname{Per}}
\newcommand{\Perp}{\operatorname{Perp}}
\newcommand{\PrePer}{\operatorname{PrePer}}
\newcommand{\PGL}{\operatorname{PGL}}
\newcommand{\Pic}{\operatorname{Pic}}
\newcommand{\Prob}{\operatorname{Prob}}
\newcommand{\Proj}{\operatorname{Proj}}
\newcommand{\Qbar}{{\bar{\QQ}}}
\newcommand{\rank}{\operatorname{rank}}
\newcommand{\Rat}{\operatorname{Rat}}
\newcommand{\Res}{{\operatorname{Res}}}
\newcommand{\Resultant}{\operatorname{Res}}
\renewcommand{\setminus}{\smallsetminus}
\newcommand{\sgn}{\operatorname{sgn}}
\newcommand{\SL}{\operatorname{SL}}
\newcommand{\Span}{\operatorname{Span}}
\newcommand{\Spec}{\operatorname{Spec}}
\renewcommand{\ss}{{\textup{ss}}}
\newcommand{\stab}{{\textup{stab}}}
\newcommand{\Stab}{\operatorname{Stab}}
\newcommand{\Support}{\operatorname{Supp}}
\newcommand{\Sym}{\operatorname{Sym}}  
\newcommand{\tors}{{\textup{tors}}}
\newcommand{\Trace}{\operatorname{Trace}}
\newcommand{\trianglebin}{\mathbin{\triangle}} 
\newcommand{\tr}{{\textup{tr}}} 
\newcommand{\UHP}{{\mathfrak{h}}}    
\newcommand{\Wander}{\operatorname{Wander}}
\newcommand{\<}{\langle}
\renewcommand{\>}{\rangle}

\newcommand{\pmodintext}[1]{~\textup{(mod}~#1\textup{)}}
\newcommand{\ds}{\displaystyle}
\newcommand{\longhookrightarrow}{\lhook\joinrel\longrightarrow}
\newcommand{\longonto}{\relbar\joinrel\twoheadrightarrow}
\newcommand{\SmallMatrix}[1]{%
  \left(\begin{smallmatrix} #1 \end{smallmatrix}\right)}
  
  \def\({\left(}
\def\){\right)}


\title
{Some counting questions for matrix products}
\author[M. Afifurrahman]{Muhammad Afifurrahman \orcidlink{0000-0002-7400-320X}}

\address{School of Mathematics and Statistics, University of New South
Wales, Sydney NSW 2052, Australia}
\email{m.afifurrahman@unsw.edu.au}
\address[]{}
\email{}

\thanks{}
\begin{abstract}Given a set $X$ of $n\times n$ matrices and a positive integer $m$, we consider the problem of estimating the cardinalities of the product sets $A_1 \ldots A_m$, where $A_i\in X$. When $X=\cM_n(\Z;H)$, the set of $n\times n$ matrices with integer elements of size at most $H$, we give several bounds on the cardinalities of the product sets. 
While proving the bounds, we also give some bounds on the cardinalities of the set of solutions of the related equations such as $A_1 \ldots A_m=C$ and $A_1 \ldots A_m=B_1 \ldots B_m$. 

We also consider the case where $X$ is the subset of matrices in $\cM_n(\F)$, where $\F$ is a field, with bounded rank $k\leq n$. In this case, we completely classify the related product set.
\end{abstract}

\maketitle

\tableofcontents

\section{Introduction} \subsection{Set-up and motivation}
For positive integers $n$ and $H$, let $\Mcal_n(\Z;H)$ be the set of $n\times n$ integer matrices $$A=(a_{ij})_{i,j=1}^n$$ with $|a_{ij}|\leq H$ for $i,j=1,\ldots,n$. We see that $\Mcal_n(\Z;H)$ is of cardinality $(2H+1)^{n^2}$.

Here, we consider some questions in \textit{arithmetic statistics} of the matrices in $\Mcal_n(\Z;H)$ (where $H\to \infty$) and $\Mcal_n(\F)$ (where $\F$ is a field), where $n$ is fixed. 
In particular, we are interested on the product set \begin{align*}
	\{A_1\ldots A_m \colon A_i \in X_i\},
\end{align*} where $X_i \subseteq \Mcal_n(\Z;H)$ or $X_i \subseteq \Mcal_n(\F)$ for $i=1,\ldots, m$.

For the first case, we derive some bounds on the cardinality of the set where $X_i=\Mcal_n(\Z;H)$, and also give some bounds on the related equations over $\Mcal_n(\Z;H)$. For the second case, we completely describe the set in the case where $X_i$ is the set of matrices with rank at most $k_i\leq n$.

We first note that several similar problems over matrix rings with additive combinatorics flavor are also studied with various direction, mostly over matrices in finite fields. For some examples, \cite{FHLOS} studies the distribution of singular and unimodular matrices over subsets of a matrix ring, and \cite{KKPSL,HH1,HH2,TV} study the expansion phenomena over matrix rings. It is worth noting that \cite{MPW,NV,XG} study the number of solutions of some equations in $\cM_n(\F_q)$, where $\F_q$ is the finite field whose cardinality is a prime power $q$.
 
  In the other direction, Ahmadi and Shparlinski \cite{AhmShp} study the arithmetic statistics of matrices in $\cM_{m,n}(H;\F_p)$, the set of $m\times n$ matrices in $\F_p$ (where $p$ is prime) whose entries are bounded by $H$. Also, El-Baz, Lee, and Str{\"o}mbergsson \cite{E-BLS} asymptotically count the number of $d\times n$ integer matrices with a bounded norm whose rank (reduced modulo a prime $p$) is a fixed number $r$.

\subsection{Integer matrices with bounded entries}
For $\Mcal_n(\Z;H)$, we count the number of tuples of  matrices in $\Mcal_n(\Z;H)$ that satisfy some given multiplicative relations. This can be seen as a dual of the results of Ostafe and Shparlinski \cite{OS}, who consider counting the numbers of matrices in $\cM_n(\Z;H)$ that are multiplicatively dependent (that is, they satisfy a certain multiplicative relation). 

Our main problem for $\cM_n(\Z;H)$ is to estimate the cardinality of \begin{align}\label{eq:wcal}
    \Wcal_{m,n}(\Z;H)=\{ A_1 \ldots A_m \colon A_1, \:  \:\ldots  ,A_m \in \cM_n(\Z;H)\}
\end{align} for fixed $m$ and $n$.

For $n=1$, the question of estimating $\#\Wcal_{m,n}(\Z;H)$ is known as \textit{Erd\H{o}s' multiplication table problem}, popularised by Erd\H{o}s \cite{E}. Ford \cite{Ford} gives an asymptotic formula for this quantity when $m=2$, and Koukoulopoulos \cite{Kou} generalises Ford's result for all $m>2$. Therefore, this set-up can be seen as a noncommutative analogue of this problem.

 We first observe that $\cM_{n}(\Z;H)\subseteq \Wcal_{m,n}(\Z;H)$. This alone gives a trivial lower bound \begin{align*}
     \# \Wcal_{m,n}(\Z;H) \gg H^{n^2}.
 \end{align*}
Moreover, by noting there are $O(H^{n^2})$ choices for each of $A_1$, $\ldots$, $A_m$, we have a trivial upper bound \begin{align*}\# \Wcal_{m,n}(\Z;H)= O(H^{mn^2})\end{align*} for all $m$. These two bounds are improved in Theorem~\ref{thm:Wmn}.

While proving Theorem~\ref{thm:Wmn}, we also give some bounds on the number of solutions of related equations in $ \cM_n(\Z;H)$, such as \begin{align}\label{eq:1}A_1 \ldots A_m=C
\end{align} for a fixed matrix $C$, and \begin{align}\label{eq:2}A_1 \ldots A_m=B_1 \ldots B_m.
\end{align}
We consider these problems over both the sets $\cM_{n}(\Z;H)$ and $\cM_{n}^*(\Z;H)$, the set of nonsingular matrices in $\cM_{n}(\Z;H)$.

When $n=1$, these problems are equivalent to bounding the number of divisors of a positive integer $k$. In this case, the result is well-known; see Lemma~\ref{lem:div}. However, the situation becomes trickier in the bigger dimension, due to matrix noncommutativity and the absence of a prime number factorisation analogue for matrices.

 We now define some related notations. For a set of $n\times n$ matrices $\cM$ and a matrix $C$, denote \begin{align}\label{def} \cT_m(\cM,C) = \{(A_1,  \ldots, \: A_m) \in \cM^m \colon A_1 \ldots A_m=C \}.
	\
\end{align} Also, define \begin{align*}
    \cT_{m}(\cM)= \{(A_1,\:\ldots,\:A_m,  B_1,\:\ldots,\:B_m)  \in \cM^{2m}: A_1\ldots A_m = B_1\ldots B_m\}.
\end{align*}

Using these, we see that the number of solutions of Equations~\eqref{eq:1} and \eqref{eq:2} are  $\#\cT_m(\cM_{n}(\Z;H),C)$ and  $\#\cT_m(\cM_{n}(\Z;H))$, respectively.
We first note that for any $C$, $$\#\cT_m(\cM_{n}(\Z;H),C) \ll H^{mn^2-1}$$ by fixing all entries of $A_2,\ldots,A_{m}$ and all, except one, entries of $A_1$. Furthermore, if $C$ is nonsingular, we see that any matrices $A_2,\ldots,A_m$ lead to at most one matrix $A_1$ that satisfies Equation~\eqref{eq:1}. Therefore, in this case, $$\#\cT_m(\cM^*_{n}(\Z;H),C)=\#\cT_m(\cM_{n}(\Z;H),C)\ll H^{(m-1)n^2}.$$ 
By using the same argument, we see that  \begin{align*}
 \#\cT_m(\Mcal^*_n(\Z;H)) \leq \#\cT_m(\Mcal_n(\Z;H)) \ll  H^{2mn^2}.
\end{align*} These trivial upper bounds are improved in Theorem~\ref{thm:21} and  Corollary \ref{thm:24}.

For the lower bounds, we observe that \begin{align*}
	\#\cT_m(\Mcal^*_n(\Z;H))\gg H^{mn^2}
\end{align*} by taking $A_i=B_i$ for $i=1,\ldots,m$ and noticing that there are asymptotically $H^{n^2}$ matrices in  $\cM^*_n(\Z;H)$. Also, \begin{align*}
	\#\cT_m(\Mcal_n(\Z;H))\gg H^{(2m-2)n^2}
\end{align*} by taking $A_m=B_m=O_n$ in the related equation.

 \subsection{Matrices over arbitrary field}
For $\cM_n(\F)$, we let $S_1,   \ldots, S_m \subseteq \Mcal_n(\F)$ be some subsets of $ \Mcal_n(\F)$ with some prescribed properties. We are interested in studying the sets of the form \begin{align*}
     \{A_1\ldots A_m \colon A_i \in S_i\text{ for }i=1,\ldots,m \}.
 \end{align*}

In this paper, we explore the question where $S_i$, with $i=1,\ldots,m$, are the sets of all $n\times n$ matrices with bounded rank $k_i\leq n$. This is done in Theorem~\ref{thm:27} and Corollary~\ref{thm:28}, using a matrix construction in Theorem~\ref{thm:26}.

\subsection{Notations}
We recall that the notations $U= O(V)$, $U\ll V$ and
$V\gg U$ are equivalent to $|U| \leq cV$ for some positive constant $c$, which
throughout this work, all implied constants may depend only on $m$ and $n$. On the other hand, we write $U(x)=o(V(x))$ if  $\lim_{x\to \infty} (U(x)/V(x))=0$. We also write $U \sim V$ if $U=O(V)$ and $V=O(U)$. We also write $U=V^{o(1)}$ if, for a given $\varepsilon>0$, we have $V^{-\varepsilon}\leq |U| \leq V^{\varepsilon}$ for large enough $V$. 

We denote $O_{n}$ as the zero $n\times n$ matrix and $\textbf{0}$ as the zero vector (whose size is taken depending on the context). 
\section{Main results} 
\subsection{Integer matrices with bounded entries}
We first consider the equation $A_1 \ldots A_m=C$, where $C $ is fixed.
We obtain the following uniform bounds for the cardinality of  $\cT_m(\cM_{n}(\Z;H),C)$ (defined in Equation~\eqref{def}).

\begin{theorem}\label{thm:21}  Let $C$ be a nonsingular matrix. Then, uniformly
     \begin{align*}
         \#\cT_m(\cM_{n}(\Z;H),C)   \leq \begin{cases}
         	H^{(m-1)(n^2-n)+o(1)},&\text{ if }C\text{ is nonsingular, for all } m,\\
         	 H^{n^2+o(1)},&\text{ if }C\neq O_n\text{ is singular and }m=2,\\
         	 H^{mn^2-n},&\text{ if }C\neq O_n\text{ is singular and }m\geq 3.
         \end{cases} 
     \end{align*}
     Furthermore, when $C=O_n$ where $O_n$ is the zero $n\times n$ matrix, we have 
   \[H^{(m-1)n^2} \ll \#\cT_m(\cM_{n}(\Z;H),O_{n}) \leq H^{(m-1)n^2+o(1)}. \]
\end{theorem}

It would be interesting to tighten these bounds. For the case $C=O_n$, it is also interesting to obtain a better error term on the asymptotical cardinality of $\#\cT_m(\cM_{n}(\Z;H),O_{n})$, as $H\to \infty$. 

Next, we give some upper bounds on the number of solutions of the equation \begin{align*}A_1\ldots A_m = B_1\ldots B_m,\end{align*} where either all matrices are in $\cM_n(\Z;H)$ or all are in $\cM^*_n(\Z;H)$ (in other words, all matrices are nonsingular). We obtain the following bounds as a corollary of Theorem~\ref{thm:21}.

\begin{corollary}    \label{thm:24}For all $m,n\geq 2$, we have
	\begin{align*}
		\#\cT_m(\Mcal^*_n(\Z;H))\leq H^{(2m-1)n^2-(m-1)n+o(1)}.
	\end{align*}We also have \begin{align*}  \#\cT_2(\Mcal_n(\Z;H))\leq  \begin{cases}
		H^{3n^2-n+o(1)},&\text{ if }m=2,\\
		H^{4mn^2-2n+o(1)},&\text{ if }m\geq 3.
	\end{cases} 
	\end{align*} 
\end{corollary}  

Indeed, it would be interesting to reduce the gap between the upper and lower bounds of these quantities.

Finally, we give some bounds on $\#\Wcal_{m,n}(\Z;H)$, defined in Equation~\eqref{eq:wcal}.
 We obtain the following result.

\begin{theorem}\label{thm:Wmn}
For all $m,n\geq 2$, we have \begin{align*}
	H^{n^2+mn-n+o(1)}\leq \#\Wcal_{m,n}(\Z;H).
	\end{align*} If $m\geq 6$, we also have \begin{align*}\#\Wcal_{m,n}(\Z,H)=o(H^{mn^2}).\end{align*}
\end{theorem}
It would be very interesting to reduce the gap between these two bounds. Furthermore, we also believe that the upper bound remains true for $2\leq m \leq 5$.

 \subsection{Matrices over an arbitrary field}
In the spirit of the problems for $\cM(\Z;H)$ in the preceding section, we first consider the set \begin{align*}
    \{ AB \colon A \in S_1,\: B\in S_2
    \},
\end{align*} where $S_1,\:S_2$ are some subsets of $\cM_n(\F)$ with some prescribed properties. In this section, we consider the sets $S_1$ and $S_2$ as the set of matrices in $\cM_n(\F)$ with some bounded rank $k\leq n$.

More precisely, let $\cM_n(\F;k)$ denote the set of matrices in $\cM_n(\F)$ with rank at most $k$. We see that $\cM_n(\F;n)=\cM_n(\F)$ and $\cM_n(\F;n-1)=\cM_n^{0}(\F)$, the set of all singular matrices in $\cM_n(\F)$. Denote \begin{align*}
S_{n}(\F;k_1,k_2)=\{
AB \colon  A\in \cM_n(\F;k_1),B\in \cM_n(\F;k_2)
\}.
\end{align*}

We prove the following assertion:

\begin{theorem}\label{thm:27} Let $n$ be a positive integer. For any nonnegative integers $k_1$ and $k_2$ with $k_1,k_2\leq n$, we have that
$$ S_{n}(\F;k_1,k_2)=  \cM_n(\F;\min{(k_1,k_2)}).$$ 
\end{theorem}
In particular, by setting $k_1=k_2=n-1$ in Theorem~\ref{thm:27}, we have that any singular matrices in $\cM_n(\F)$ can be represented as a product of two singular matrices.

By inducting on $m$, Theorem~\ref{thm:27} can be generalised as follows:
\begin{corollary}\label{thm:28} For any nonnegative integers $m,\:n\geq 2$ and $ k_1, \ldots,\: k_m \leq n$, we have that
$$ \{
A_1 \ldots A_m \colon  A_i\in \cM_n(\F;k_i)
\}= \cM_n(\F;\min{(k_1, \:\ldots,\:k_m)}). $$
\end{corollary}

We end this subsection with noting that if $\F$ is a finite field, then the formula for $\#\cM_n(\F;k)$ is known; see Fisher~\cite{Fisher}.

\section{Preliminary results}
We first need a result on the equivalence of counting $\cW_{m,n}(\Z;H)$  for $n=1$; namely, counting the cardinality of the set \begin{align*}
  \cA_m(H)=  \{a_1 \ldots a_n \colon a_1, \ldots, a_m \in \Z, |a_i|\leq H\text{ for }i=1,\ldots,m\}.
\end{align*}
By a result of Koukoulopoulos \cite[Corollary~1]{Kou} on the product sets of positive integers not less than $H$, we have the following result.
\begin{lemma}\label{thm:F}For $m\geq 2$, let $\rho=(m+1)^{1/m}$ and $$Q(u):=\int_1^u\log t\,dt=u\log u-u+1\quad(u>0).$$ We have that $$
\# \cA_{m}(H)\sim \frac{H^m}{(\log H)^{Q(\frac1{\log\rho})}(\log\log H)^{\frac{3}{2}}}.
$$ 
\end{lemma}

In addition, some of our arguments on bounding $\cM_n(\Z;H)$ are based on the properties of integral matrices with bounded norm. We quote the following results from Shparlinski \cite{IES} and Katznelson \cite{K}

\begin{lemma}\label{thm:KDRS}
	Fix an integer $d$. Then, uniformly there are $O(H^{n^2-n}\log H)$ matrices in $\Mcal_n(\Z;H)$ of determinant $d$. Also, when $d=0$, there are asymptotically $H^{n^2-n}\log H$ matrices in $\Mcal_n(\Z;H)$ of determinant $0$.
\end{lemma}

We note that an asymptotic formula for the quantities in Lemma~\ref{thm:KDRS} exists from Duke, Rudnick, and Sarnak \cite[Example 1.6]{DRS}, when $d\neq 0$ is fixed. However, this bound cannot be used in our arguments since their argument requires $d$ to be fixed.

We also use the following property of integral matrices from Katznelson \cite[Theorem~1~(2)]{K2}.
\begin{lemma}\label{lem:33} Let $k$ be an integer with $1\leq k < n$. Then, there are $H^{nk+o(1)}$ matrices in $\Mcal_n(\Z;H)$ of rank $k$.
\end{lemma}


We also use a well-known result concerning the number of divisors $\tau(k)$ of a positive integer $k$, see for example \cite[Equation~(1.81)]{IK}.

\begin{lemma}\label{lem:div} We have that $$\tau(k)=k^{o(1)}$$ as $k\to \infty$. \end{lemma}

Finally, in order to prove Theorem~\ref{thm:27}, we need the following result on decomposition of a matrix with a fixed rank $k$. The following result is new, as per the author's knowledge.

\begin{theorem}\label{thm:26}
Let $A\in \cM_n(\F)$ with $\rank A=k$. Then, there exists $B\in \cM_n(\F)$ with $\rank B = k$ such that $BA=A$.
\end{theorem}

\section{Proof of Theorem~\ref{thm:21}}
\subsection{The case of nonsingular $C$}

We first notice that in the equation $A_1 \ldots A_m=C$, where $C$ is nonsingular, all of $ A_1, \ldots, \: A_{m} $ are also nonsingular. This implies that each choice of $ A_1, \ldots,\: A_{m-1} $ gives at most a unique matrix $A_m$. Therefore, we now only count $A_i$ for $i=1,\ldots,\: m-1$.

We know that $\det A_i\mid \det C$. Therefore, from Lemma~\ref{lem:div}, there are $$\lvert \det C \rvert 
^{o(1)}\leq H^{o(1)}$$ possible values $d$ of $\det A_i$. From Lemma~\ref{thm:KDRS}, there are at most $H^{n^2-n+o(1)}$ matrices in $\Mcal_n(\Z;H)$ with determinant $d$. Therefore, there are at most $H^{n^2-n+o(1)}$ possible matrices $A_i$, for $i=1,\ldots,\: m-1$, which satisfy the equation. Multiplying these bounds proves the initial statement.

\subsection{The case of singular nonzero $C$}
We first prove the bound for general $m$. We know that in the equation \begin{align*}A_1\ldots A_m=C,
\end{align*} at least one of $A_1$, \ldots, $A_m$ (suppose $A_1$) is singular.  By Lemma~\ref{thm:KDRS}, there are at most $H^{n(n-1)+o(1)}$ possible choices of $A_1$. Since there are $H^n$ choices for other matrices, we can multiply the bounds to complete the proof. 

We now consider the case $m=2$. We bound $\#\cW_{2,n}(H,C)$, where $C$ is a fixed nonzero singular matrix of size $n\times n$, with $n\geq 2$. Consider the equation $AB=C$. Without loss of generality, let $A$ be singular, and $\rank A=k\geq 1$.  We write $A$, $B$, and $C$ in the following form: \begin{align*}
	A=\begin{pmatrix}
		X_1
		& \rvline & V_1 \\
		\hline
		W_1 & \rvline &
		Y_1
	\end{pmatrix},\: 	B=\begin{pmatrix}
		X_2
		& \rvline & V_2 \\
		\hline
		W_2 & \rvline &
		Y_2
	\end{pmatrix},\: 
C=\begin{pmatrix}
		X
		& \rvline & V \\
		\hline
		W & \rvline &
		Y
	\end{pmatrix},
\end{align*} where $X$, $X_1$, and $X_2$ 
 are $k\times k$ matrices, $Y$, $Y_1$, and $Y_2$  are $(n-k)\times (n-k)$ matrices, $V_1$, $V_2$, and $V$ are $k\times (n-k)$ matrices, and $W_1$, $W_2$, and $W$ are $(n-k)\times k$ matrices. Since $\rank A=k$, we know that $A$ has at least a nonsingular $k\times k$ submatrix. Without loss of generality, we may assume that $X_1$ is nonsingular. We then consider the following matrix equation: \begin{align*}
AB=C \iff \begin{pmatrix}
		X_1X_2+V_1W_2
		& \rvline & X_1V_2+V_1Y_2 \\
		\hline
		W_1X_2+Y_1W_2 & \rvline &
		W_1V_2+Y_1Y_2
	\end{pmatrix}=\begin{pmatrix}
		X
		& \rvline & V \\
		\hline
		W & \rvline &
		Y
	\end{pmatrix}.
\end{align*}

In particular, looking at the first row of the above equations, we obtain the following system of equations.
\begin{align*}
X_2&=X_1^{-1}(X-V_1W_2),\\
V_2&=X_1^{-1}(V-V_1Y_2).
\end{align*}

The first equation implies that fixed $A$, $C$, and $W_2$ give a unique $X_2$. The second equation implies that fixed $A$, $C$, and $Y_2$ give a unique $V_2$. Therefore,  for a fixed $A$, $W_2$, and $Y_2$, we get at most one possible matrix $B$ that satisfy the equation $AB=C$.

By Lemma~\ref{lem:33}, there are $H^{nk+o(1)}$ matrices $A\in \cM_n(\Z;H)$ such that $\rank A=k$. On the other hand, there are $O(H^{(n-k)k})$ possible choices of $W_2$ and $O(H^{(n-k)^2})$ possible choices of $Y_2$.

These observations implies that for a fixed $C$, there are at most $$H^{nk+o(1)} \cdot H^{(n-k)k} \cdot H^{(n-k)^2}  = H^{n^2+o(1)}$$ different pairs $(A,B)$ such that $AB=C$ and $\rank A=k$. Summing all possible $k$ completes the proof.

\subsection{The case of $C=O_n$}
The lower bound can be attained by considering that any tuples of matrices $(A_1,\ldots,A_m)$ with $A_1=O_n$, where $O_n$ is the zero $n\times n$ matrix, satisfies the equation \begin{align*}
	A_1\ldots A_m=O_n.
\end{align*} This implies that there are at least $H^{(m-1)n^2}$ different possible tuples of matrices that satisfy the equation above.

Now, we prove the upper bound. Suppose that $(A_1,\ldots,A_m)$ satisfy the equation above. From Sylvester's rank inequality and induction, we have \begin{align*}
\rank A_1 \ldots A_m\geq \sum_{i=1}^{n} \rank A_i - (m-1)n.
\end{align*}

In our case, since $\rank O_n=0$, this inequality implies \begin{align*}
	(m-1)n\geq \sum_{i=1}^{m} \rank A_i.
\end{align*}

Now, we let $\rank A_i=k_i$ for $i=1,\ldots, m$. We recall that from Lemma~\ref{lem:33}, there are $H^{nk+o(1)} $ matrices in $\cM_n(\Z;H)$ of rank $k$. Therefore, from the inequality above, \begin{align*}
\#\cT_m(\cM_{n}(\Z;H),O_{n}) &\leq \sum_{\substack{0\leq k_1, \ldots,k_m \leq n\\
		k_1+ \ldots+k_m\leq (m-1)n}} H^{nk_1+o(1)}\cdot \ldots \cdot  H^{nk_m+o(1)}\\
&\leq \sum_{k_1+ \ldots+k_m\leq (m-1)n} H^{n(k_1+\ldots+k_m)+o(1)} \\
&\leq \sum_{k_1+ \ldots+k_m\leq (m-1)n} H^{n(m-1)n+o(1)} \\
&= H^{(m-1)n^2+o(1)} .
\end{align*} This completes the proof of the upper bound.

\section{Proof of Corollary~\ref{thm:24}}
\subsection{The nonsingular case}

Consider the equation \begin{align}\label{eq:AB}
    A_1 \ldots A_m=B_1 \ldots B_m,
\end{align} where $A_i, B_i \in \cM_n^*(\Z;H)$ for $i=1,\ldots,m$. For each choice of $( B_1,   \ldots, B_m )$, we see that, by Theorem~\ref{thm:21}, there are at most $O(H^{(m-1)(n^2-n)})$ possible $m$-tuples of matrices $(A_1,\:\ldots,\:A_m)$ that satisfy Equation~\eqref{eq:AB}. On the other hand, there are at most $O(H^{n^2})$ possible matrices $B_i$, for $i=1,\ldots,\: m$. Multiplying these bounds proves the statement.

\subsection{The general case}We first give an upper bound for the case $m\geq 3$. Notice that, in contrast to the previous case, we now allow the matrices to be singular. We first consider the case where all matrices in the equation \begin{align*}
A_1\ldots A_m=B_1\ldots B_m
\end{align*} are nonsingular. Recalling the result in the previous section, there are at most $H^{(2m-1)n^2-(m-1)n+o(1)}$ solutions in this case.

Next, we consider the case where at least one of them is singular. Without loss of generality, suppose $A_1$ is singular. This implies that at least one of $B_1$, $\ldots$, $B_m$ (without loss of generality, $B_1$) is also singular. Then, from Lemma~\ref{lem:33} there are at most $H^{n(n-1)+o(1)}$ choices for $A_1$ and $B_1$. Since there are $H^n$ possible choices for other matrices, there are $H^{4mn^2-2n+o(1)}$ possible solutions in this case. Combining both cases, we see that 
\begin{align*}
	\#\cT_n(\Mcal_n(\Z;H))\ll \max\{H^{4mn^2-2n+o(1)}, H^{(2m-1)n^2-(m-1)n+o(1)}\} =H^{4mn^2-2n+o(1)}.\end{align*}
	
We now proceed to improve this bound for $m=2$. 
Using the argument for general $m$, we see that  there are at most $H^{3n^2-n+o(1)}$ solutions to the equation $AB=CD$ where all of them are nonsingular. Next, we count the solutions of $AB=CD$, where at least one of $A$, $B$, $C$, and $D$ is singular. Without loss of generality, we fix $CD$, where $C$ is singular. By Lemma~\ref{thm:KDRS}, this can be achieved in $H^{2n^2-n+o(1)}$ ways. By Theorem~\ref{thm:21}, we see that in this case, there are at most $H^{n^2+o(1)}$ possible pairs of matrices $(A,B)$ such that $AB=CD$. Multiplying these bounds, we find that there exist at most $H^{3n^2-n+o(1)}$ solutions of the equation $
AB=CD$, where at least one of $A$, $B$, $C$, and $D$ is singular. 

Combining both cases, we see that 
\begin{align*}
  \#\cT_2(\Mcal_n(\Z;H))\ll H^{3n^2-n+o(1)},
	\end{align*} which completes the proof.

\section{Proof of Theorem~\ref{thm:Wmn}}
\subsection{Upper bound}
We will prove a stronger bound than what is stated in Theorem~\ref{thm:Wmn}, namely \begin{align}\label{eqn:kou}
	\#\Wcal_{m,n}(\Z;H) = O\left(\frac{H^{mn^2}}{(\log H)^{Q(\frac1{\log\rho})-1}(\log\log H)^{\frac{3}{2}}}\right),
\end{align} where $\rho=m^{1/(m-1)}$ and 
	$Q(u):u\log u-u+1.$

We first notice that the absolute values of the entries in $\Wcal_{m,n}(\Z;H)$ are bounded above by $n^{m-1}H^m$. Hence, $\Wcal_{m,n}(\Z;H)\subseteq \Mcal_{m,n}(\Z;n^{m-1}H^m)$. 

Now, let $\cD$ be the set of all possible values of $\det A_1\ldots \det A_m$.
We recall that $|\det A|\leq n^{m-1}H^m$ for any $A \in \cM_n(\Z;H)$.  Therefore, by Lemma~\ref{thm:F}, \begin{align*} \label{eq:cd}	\#\cD &= O\left(\frac{(H^n)^m}{(\log H^n)^{Q(\frac1{\log\rho})}(\log\log H^n)^{\frac{3}{2}}}\right)\\& = O\left(\frac{H^{mn}}{(\log H)^{Q(\frac1{\log\rho})}(\log\log H)^{\frac{3}{2}}}\right).
\end{align*}

From Lemma~\ref{thm:KDRS},  there are $O(H^{m(n^2-n)}\log H)$ matrices in $\Mcal_{m,n}(\Z;n^{m-1}H^m)$ with determinant $d$, for each $d\in \cD$. Therefore, we have \begin{align*}
	\cW_{m,n}(\Z;H)&\leq \#\cD \cdot O(H^{n^2-n}\log H)\\
	&\leq O\left(\frac{H^{mn^2}}{(\log H)^{Q(\frac1{\log\rho})-1}(\log\log H)^{\frac{3}{2}}}\right),
\end{align*}
which completes the proof of Equation~\eqref{eqn:kou}.

It remains to prove that Equation~\eqref{eqn:kou} implies $\#\Wcal_{m,n}(\Z;H)=o(H^{mn^2})$ for $m\geq 6$. To prove this, we notice that the function $f(m)=Q\left(\frac{1}{\log (m^{1/(m-1)})}\right)$ is increasing on $m$, and $f(6)\geq 1$. Therefore, in this case, the denominator of the right-hand side of Equation~\eqref{eqn:kou} would imply $$\#\Wcal_{m,n}(\Z;H)=o(H^{mn^2}),$$ for $m\geq 6$, proving the desired upper bound.

We end this proof with noting that  Equation~\eqref{eqn:kou} remains true for $2\leq m \leq 5$. However, since $f(5)$ is smaller than $1$, that equation alone does not imply $\#\Wcal_{m,n}(\Z;H)=o(H^{mn^2})$ in this case.

\subsection{Lower bound}We consider the lower bound of $\#\Wcal_{m,n}^*(\Z;H)$. From Theorem~\ref{thm:21}, there are at most $H^{(m-1)(n^2-n)+o(1)}$ different possible solutions to the equation $A_1 \ldots A_m=C$, where $A_i \in \cM_n(\Z;H)$ (for $1\leq i \leq m$) and $C$ is a fixed nonsingular matrix.

We also see that there are asymptotically $H^{n^2}$ nonsingular matrices in $\cM_n(\Z;H)$. Therefore, if we count the number of $(m+1)$ -tuples $(A_1,\ldots,A_m,C)$ that satisfies $A_1 \ldots A_m=C$ and $C\in \Wcal^*_{m,n}(\Z;H)$, we get that \begin{align*}
	(\#\Wcal^*_{m,n}(\Z;H))H^{(m-1)(n^2-n)+o(1)} &\geq H^{mn^2+o(1)} \\  \iff   \#\Wcal^*_{m,n}(\Z;H) &\geq 
	H^{n^2+mn-n+o(1)}.
\end{align*} Since $\Wcal^*_{m,n}(\Z;H)\subseteq \Wcal_{m,n}(\Z;H)$, the lower bound is proven.

\section{Proof of Theorem~\ref{thm:27}}

Let $n\geq 2$ be a positive integer and $k_1,k_2$ be nonnegative integers with $0\leq k_1,k_2\leq n$. Without loss of generality, suppose $k_1 \leq k_2$. Since \begin{align*}
    \rank AB\leq \min{(\rank A, \rank B)}
\end{align*} for any matrices $A$ and $B$, we see that the rank of a matrix in $S_n(\F;k_1,k_2)$ is at most $k_1$. This implies that $S_n(\F;k_1,k_2)\subseteq \cM_n(\F;k_1)$.

Now we prove that any matrix $C$ in $\cM_n(\F;k_1)$ can be represented as a product of two matrices $XY$, where $\rank X \leq k_1$ and $\rank Y \leq k_2$. However, by Theorem~\ref{thm:26}, there exists a matrix $B$ with $\rank B= k_1\leq k_2$ that satisfies \begin{align*}
    BC=C.
\end{align*}

This implies that $\cM_n(\F;k_1)\subseteq S_n(\F;k_1,k_2)$, which concludes the proof.

\section{Proof of Theorem~\ref{thm:26}}
If $k=n$, we pick $B=I$. If $k=0$, we pick $B=O_n$. Suppose that $1\leq k<n$. Since $\rank A=k$, $\text{Null}(A)=n-k$. This implies that there exist $n-k$ linearly independent vectors $\mathbf{v}_1, \ldots,\:\mathbf{v}_{n-k} \in \F^n$ such that $\textbf{v}_i^\intercal A=\textbf{0}$ for all $i=1,\:\ldots,\:n-k$.

Consider an $ (n-k)\times n$ matrix $C'$ with \begin{align*}
    C'=[\textbf{v}_1 \mid  \ldots \mid \textbf{v}_{n-k}]^\intercal.
\end{align*}
Let its reduced row echelon form be \begin{align*}
    C=[\textbf{w}_1 \mid  \ldots \mid \textbf{w}_{n-k}]^\intercal.
\end{align*}

First, we see that the vectors $\{\textbf{w}_1,  \ldots , \textbf{w}_{n-k}\}$ are linearly independent. Also, $\textbf{w}_i^\intercal A=\textbf{0}$ for $i=1,\ldots, n-k$.

On the other hand, we observe the following properties of $C$ from the definition of the reduced row echelon form: \begin{enumerate}
    \item The leading nonzero components of $\textbf{w}_i$ is 1. Suppose that the leading 1 of $\textbf{w}_i$ in $C$ is in column $z_i$.
    \item If a column of $C$ has a leading 1, its other components are zero, and
    \item $z_1<\ldots<z_{n-k}$.
\end{enumerate}

Now define $\mathbf{b}_{z_i}=-\mathbf{w}_i$ for $i=1,\ldots,\:n-k$, with $z_i$ defined as above. Also, let $\mathbf{b}_j=\textbf{0}$ for other index $j$ such that $\textbf{b}_j$ is defined for all $j=1,\:\ldots,\:n$.

Let $$B'=[\mathbf{b}_1 \mid \ldots \mid \mathbf{b}_n]^\intercal.$$ We notice that $B'_{z_iz_i}=-1$ for $i=1,\ldots,\:n-k$. Also, $B'_{jz_i}=0$ if $j\neq z_i$. Furthermore, since $\textbf{b}_{i}^\intercal A=\textbf{0}$ for $i=1,\:\ldots,\:n$, we have $B'A=O_n$. Therefore, letting $B=B'+I$, we have $$ BA = (B' + I) A = A.$$ This implies that $B$ satisfies the first part of our theorem.

Now we calculate the rank of $B$. By the definition of the vectors $\mathbf{b}_i$, we see that the $z_i$-th column of $B$ is $\textbf{0}$, for $i=1,\:\ldots,\:n-k$. Then, $B$ has at least $n-k$ zero column and $\rank B\leq k$. On the other hand, if $j\neq z_i$ for any $i$, the $j$-th row of $B$ consists only of a $1$ in its $j$-th component and zeros in other components. Therefore, these $k$ rows are linearly independent, which implies $\rank B \geq k$. Therefore, $\rank B=k$, which concludes the proof.
\section*{Acknowledgements} The author would like to thank Alina Ostafe and Igor Shparlinski for many ideas, comments, and corrections during the preparation of this work. The author would also like to thank anonymous reviewers for their suggestions. The author is supported by UNSW Tuition Fee Scholarship and Australian Research Council Grant DP200100355.

\end{document}